\documentclass{article}
\usepackage[english]{babel}

\usepackage{amssymb,amsmath,graphics,tikz}

\usepackage{hyperref}
\usepackage{enumerate
}
\usepackage{todonotes}

\parskip\medskipamount
\newtheorem{theorem}{Theorem}[section]

\newtheorem{remark}[theorem]{Remark}
\newtheorem{prop}[theorem]{Proposition}
\newtheorem{lemma}[theorem]{Lemma}

\def\<{\langle}
\def\>{\rangle}
\newcommand{\proof}{\emph{Proof.~}}

\def\qed{{\hfill\hphantom{.}\nobreak\hfill$\Box$}}

\newcommand{\R}{\mathbb{R}} 
\newcommand{\N}{\mathbb{N}}

\newcommand{\C}{\mathbb{C}}
\newcommand{\HH}{\mathbb{H}}
\newcommand{\OO}{\mathbb{O}}

\newcommand{\A}{\mathbb{A}} 
\newcommand{\B}{\mathbb{B}}

\newcommand{\SO}{\mathsf{SO}}
\newcommand{\SU}{\mathsf{SU}}
\newcommand{\Aut}{\mathrm{Aut}}

\DeclareMathOperator{\Pu}{Pu}
\let\Re\relax
\DeclareMathOperator{\Re}{Re}

\usepackage{mathtools}

\DeclarePairedDelimiter\floor{\lfloor}{\rfloor}

\newcommand{\define}{\mathrel{\mathop:}=}

\begin{document}

\title{On exceptional homogeneous compact geometries of type $\mathsf{C}_3$}
\author{Jeroen Schillewaert \\  \small{jschillewaert@gmail.com} \and Koen Struyve \\ \small{kstruy@gmail.com}}
\date{ \today }

\maketitle

\begin{abstract}
We provide a uniform framework to study the exceptional homogeneous compact geometries of type $\mathsf{C}_3$. 
This framework is then used to show that these are simply connected, answering a question by Kramer and Lytchak, and to calculate the full automorphism groups.
\end{abstract}

\small{
}
\section{Introduction}

In~\cite{Tit:81} Tits introduced geometries of Coxeter type, which are geometries which locally are buildings. 
In particular, buildings are examples of such geometries. The motivation for these objects comes from applications to (finite) group theory, see for example~\cite[C.7]{Kan:86}.

Flag-transitive finite geometries of type $\mathsf{C}_3$ have been classified by Aschbacher and Yoshiara in~\cite{Asc:84} and~\cite{Yos:96}. 
Such a geometry is either a building, or it is isomorphic with the Neumaier geometry on seven points (see~\cite{Neu:84}).

In the compact connected case a similar classification has been obtained by Kramer and Lytchak in~\cite{Kra-Lyt:14} (albeit missing one of the two exceptional cases). 
The eventual conclusion here is that such a geometry is either covered by a building, or it is isomorphic to one of two exceptional geometries. 
It are precisely these geometries which are the subject of the current paper. 

The two exceptional geometries were first encountered in the study and classification of polar actions on manifolds, see~\cite{Pod-Tho:99},~\cite{Fan-Gro-Tho1}, ~\cite{Fan-Gro-Tho2} and~\cite{Gor-Kol}. They are associated to actions of $\SU(3) \times \SU(3)$ and $\SO(3) \times \mathsf{G}_2$ on the Cayley plane $\mathbb{OP}^2$.

In this paper we provide a uniform description for both geometries, and use this to show that these geometries are simply connected. This answers Problem 5 of~\cite{Kra-Lyt:14}. We also use this description to obtain their full automorphism groups.

In Sections~\ref{section:comp} and~\ref{section:geom} we introduce basic notions needed later on. In Section~\ref{section:geomcomp} both of the exceptional geometries are constructed in an uniform way, and we discuss briefly their structure. 
Section~\ref{section:auto} discusses the automorphism group, the last section~\ref{section:simple} is devoted to the proof of their simple connectedness.

{\bf Acknowledgement.} The first author wants to thank Linus Kramer profoundly for the hospitality and the interesting mathematical discussions during his stay at the University of M\"unster as a Von Humboldt Fellow. Both authors would like to thank Alexander Lytchak for helpful comments and discussions on this problem.

\section{Composition algebras}\label{section:comp}

Let $\A$ be a composition algebra over the real numbers $\R$, so $\A$ is either $\R$ itself or $\C$, $\HH$ or $\OO$ (by Hurwitz's theorem). The algebra $\A$ comes equipped with a norm $ | \cdot |: \A \to \R$ and  standard involution $\bar{\cdot}$.  Let $k$ be either $\R$ of $\C$ and assume that $k$ is a subfield of $\A$.

The norm on $\A$ induces a positive-definite $\R$-bilinear form $\< \cdot , \cdot \> : \A \times \A \to \R$ on $\A$ when interpreted as an $\R$-vector space. Therefore we can consider orthogonality in $\A$. The orthogonal complement of a subspace $K$ in $\A$ is denoted by $\perp_\A K$.

We define  the \emph{$k$-pure} elements in $\A$ as $\Pu_k(\A)  \define \perp_\A k$. In particular we have a direct sum decomposition 

$$
\A = k \oplus \Pu_k (\A).
$$ 

Note that the standard involution acts on $\Pu_k(\A)$ as the map $x\mapsto -x$..

We now define an Hermitian inner product $( \cdot | \cdot ): \A \times \A \to k$ on $\A$ (as a right vector space) over $k$ by setting $(x | y)$ to be the $k$-part of the product $xy$. Note that if $k=\R$, then $\< \cdot , \cdot \>$ and $( \cdot | \cdot )$ agree. This Hermitian inner product has the following properties.

\begin{align}
\Re (a | b) &= \< a, b \> \\
| (a | b)| &\leq |a | \cdot |b| \\
| (a | b)| &= |a | \cdot |b| \iff a,b \mbox{ are } k\mbox{-linearly dependent} \label{eq:CS} \\
(a |a) &= |a|^2
\end{align}

where $a,b \in \A$. Proofs of these are either direct or are basic properties of inner products. For an element $a \in \Pu_k(\A)$ we denote the orthogonal complement of the vector line $ak$ in $\Pu_k(\A)$ by $a^\perp$.

The following result will be crucial for our construction.

\begin{prop}\label{prop:map}
Let $\A$ and $\B$ be two composition algebras over $\R$, both containing a common subfield $k$ which is either $\R$ or $\C$. Assume additionally that the dimension of $\A$ over $k$ is four, and where the $k$-dimension of $\B$ is at least the $k$-dimension of $\A$.
Let $a,c \in \Pu_k(\A)$ and $b,d \in \Pu_k(\B)$ be of norm 1 such that $(a | c)=(b | d)$ and $ak \neq ck$. Then there exists a $k$-algebra morphism of $\A$ into $\B$ mapping the pair $(a,c)$ to $(b,d)$.
\end{prop}
\proof
Note that $bk  \neq dk$ by Equation~\ref{eq:CS}. We may assume without loss of generality that $(a|c) = (b|d)= 0$ by replacing $c$ with a suitable scalar multiple of $c - (a | c) a$ being of unit norm (and similarly for $d$). In particular this implies that $ac$ and $bd$ are again in $\Pu_k(\A)$, see e.g. Proposition 11.10 of \cite{CPP}.

The elements $1,a,c$ and $ac$ then form a $k$-basis for $\A$, satisfying the following identities:
\begin{align*}
\overline{a} = a^{-1} &= -a \\
\overline{c} = c^{-1} &= -c \\
\overline{ac} =  \overline{c} \overline{a} &= ca = -ac \\
le = - \overline{le} &= - \overline{e} \overline{l} = e \overline{l}
\end{align*}
where $l \in k$ and $e \in \Pu_k(\A)$. 
Note that these relations completely determine $\A$ as a $k$-algebra. As analogous properties hold for $b$ and $d$, we can extend the map $(1,a,c,ac) \mapsto (1,b,d,bd)$ to a $k$-algebra morphism from $\A$ into $\B$. \qed

\section{Geometries}\label{section:geom}

A \emph{geometry over a set $I$} is a system $\Gamma \define (V,\tau, *)$, consisting of a set $V$, a surjective map $\tau: V \to I$, and a binary symmetric relation $*$ on $V$ such that for any two elements $x,y$ of $V$ whose images under $\tau$ are identical, the relation $x * y $ holds if and only if $x=y$. The relation `$*$' is the \emph{incidence relation}, the image by $\tau$ of an element or a subset of $V$ is its \emph{type}.

A \emph{flag} of the geometry $\Gamma$ is a set of pairwise incident elements of $V$. The \emph{type} of a flag is its image under $\tau$.
The \emph{corank} of a flag is $| I |$ minus the size of its type.
Let $X$ be a flag, and let $Y$ be the set of all elements of $V \setminus X$ incident to $X$. 
The system $\Gamma_X := (Y, \tau\vert_Y, * \cap (Y \times Y))$ forms a geometry over $I \setminus \tau(X)$ and is called the \emph{residue of $X$} in $\Gamma$. For a $J \subseteq I$, the set of all flags of the type $J$ as subset of $V^J$ is called the \emph{flag variety} $V_J$.

The geometry $\Gamma$ is \emph{connected} if the graph with vertices $V$ and adjacency relation `$*$' is connected. A geometry is \emph{thick} if every flag of corank 1 is contained in at least three maximal flags. 

A \emph{generalized $n$-gon} (with $2 \leq n < \infty$) is a thick rank 2 geometry such that the associated graph has girth $2n$ (i.e. the smallest cycle has length $2n$) and diameter $n$.
Generalized triangles correspond exactly with (axiomatic) projective planes.

We call the geometry $\Gamma$ a \emph{compact geometry} if the set of vertices $V$ carries a compact Hausdorff topology such that for every $J \subseteq I$ the flag variety $V_J$ is closed in the compact product space $V^J$.

\subsection{Geometries of type $\mathsf{C}_3$}
A connected geometry over the set $\{1,2,3\}$ is a \emph{geometry of type $\mathsf{C}_3$} if 
\begin{itemize}
\item
the residue of a vertex of type 1 is a generalized quadrangle,
\item
the residue of a vertex of type 2 is a generalized digon,
\item
the residue of a vertex of type 3 is a generalized triangle.
\end{itemize}

One often calls the vertices of type 1 the \emph{points}, those of type 2 \emph{lines}, and of type 3 \emph{planes}. In this sense we can speak of \emph{collinearity} of points and \emph{coplanarity} of lines. 

By a result of Tits (\cite[6.2.3]{Tit:81}), the buildings of type $\mathsf{C}_3$ are precisely those geometries of type $\mathsf{C}_3$ in which no pair of points is incident with more than one common line.

The next lemma captures part of the structure of these geometries, which is well known, see e.g.~\cite[Exercise 7.7]{Pas:94}.

\begin{lemma}\label{lem:c3}
Let $\Gamma$ be a geometry of type $\mathsf{C}_3$. 
\begin{enumerate}[i.]
\item
For a point $p$ and a line $L$ there exists at least one point $q$ incident with $L$ and collinear with $p$.
\item
For a point $p$ and plane $\pi$, there exists at least one line $L$ incident with $\pi$ and having a common incident plane with $p$.
\end{enumerate}
\end{lemma}
\proof
Since the geometry is connected and since the residue of a vertex of type 3 is a generalized triangle (i.e. a projective plane) there exists a (finite) path from $p$ to $L$ in the incidence graph
alternating between points and lines. Let $a*D*b*E*c*F$ be a path where $a,b,c$ are points and $D,E,F$ are lines. Consider a plane $\pi* F$. Since the residue of $c$ is a generalized quadrangle there exist a plane $\pi' * E$ with $\pi' \cap \pi=L_c$. Since the residue in $b$ is a generalized quadrangle there exists a plane $\pi'' * D$ with $\pi'\cap \pi"=L_b$. Then  $d\in L_b\cap L_c$ is a point in $\pi$ collinear with $a$. Considering the residue in $d$ yields there exists a plane containing $a$ intersecting $\pi$ in a line $L_d$.
Since $\pi$ is a projective plane the intersection $L_d\cap F$ contains a point $e$. Hence we have found a shorter path from $a$ to $F$, namely $a * \langle a,e\rangle * e *F$, which implies (i).

For (ii) consider any line $L$ in $\pi$. Then by (i) there is a point $q$ incident with $L$ and collinear with $p$. Considering the residue of $q$ proves the statement.
\qed

\subsection{Coverings}

Let $\Gamma\define (V,\tau, *)$ be a geometry of type $\mathsf{C}_3$. From $\Gamma$ one can construct a connected simplicial complex with vertices the elements of $V$ and simplices the flags of $\Gamma$. We denote the metric realization of this simplicial complex by $| \Gamma |$. Connected covers of $| \Gamma |$ then correspond with 2-covers of the geometry $\Gamma$ (see~\cite[Ch. 12]{Pas:94}), which are again geometries of type $\mathsf{C}_3$.

By this correspondence one can consider the \emph{universal cover} $\widetilde{\Gamma}$ of $\Gamma$, and say that the geometry is simply connected if it is its own universal cover.

\section{Geometries from composition algebras}\label{section:geomcomp}

In this section we construct a class of geometries of type $\mathsf{C}_3$ starting from two composition algebras with a common subfield. The construction itself is contained in Section~\ref{sec:constr}, the verification of the claim in Section~\ref{sec:verif}. 

Lastly, in Section~\ref{sec:pos}, we list the different possibilities for the composition algebras and show the existence of a compact flag-transitive automorphism group. 

\subsection{Construction of $\Gamma$}\label{sec:constr}

We define a geometry $\Gamma$ as follows. Let $\A$ and $\B$ be two composition algebras over $\R$, containing a common subfield $k$ which we assume to be either $\R$ or $\C$, such that $\A$ is four-dimensional over $k$ and where the ($k$-)dimension of $\B$ is at least the ($k$-)dimension of $\A$.

Scalar multiplication, as well as projectivization, is always understood to be over $k$.

\begin{itemize}
\item
The points are the vector lines in $\Pu_k(\A)$ as vector space over $k$. 
\item 
Lines are of the form $[a,b]$ where $a  \in \Pu_k(\A)$, $b \in \Pu_k(\B)$ and $|a|= |b| \neq 0$,  up to multiplication by a common scalar multiple.
\item 
The planes of the geometry are formed by the embeddings $\phi: \A \to \B$ as composition algebras over $k$.
\end{itemize}

Incidence is defined as follows. Every point is incident with every plane, a point $uk = \<u\>$ ($u \in \Pu_k(\A)$) is incident with a line $[a,b]$ if and only if $ a \in u^\perp$. A line $[a,b]$ and a plane $\phi$ are incident if and only if $\phi(a) = b$.

\begin{remark}\label{rem:flat}
Often one calls a geometry of type $\mathsf{C}_3$ where every point is incident with every plane \emph{flat} (see for example~\cite[4.3.2]{Pas:94}).
A second property that our geometry $\Gamma$ has is that if two lines are incident with at least two common points, then every point incident with one line is also incident with the other. (Another way to say this is that the point shadows of both lines agree in this case.) The set of points of $\Gamma$ and the set of point shadows of lines (with duplicity removed) form a projective plane over $k$.
\end{remark}

\subsection{$\Gamma$ is a geometry of type $\mathsf{C}_3$}\label{sec:verif}

Before we prove that $\Gamma$ is indeed a geometry of type $\mathsf{C}_3$, we obtain a criterion for coplanarity of lines. 
(Note that this criterion is not influenced by taking scalar multiples.)

\begin{lemma}\label{lem:coplanar}
Two lines $[a,b]$ and $[c,d]$ are coplanar if and only if $( a | c ) = (b | d)$.
\end{lemma}
\proof
The `only if' direction is clear as embeddings of composition algebras preserve the inner product.
The other direction follows from Proposition~\ref{prop:map}.
\qed

\begin{prop}
The geometry $\Gamma$ is of type $\mathsf{C}_3$.
\end{prop} 

\proof
Clearly $\Gamma$ is thick, connected, with the residue of any given line being a digon, and the residue of any plane being a projective plane over $k$. So it remains to verify that the residue of a given point $\< u \>$  is a generalized quadrangle.  

All planes belong to this residue and a line $[a,b]$ belongs to it if and only if $a \in u^\perp$.

In order to verify that the residue is a generalized quadrangle consider a plane $\phi$ and a line $[a,b]$ with $\phi(a) \neq b$. 
It suffices to show existence and uniqueness of a line $[c,\phi(c)]$ (so incident with $\phi$) such that $( a | c )=(b |\phi(c) )$ (up to taking a scalar product of $c$) which states that $[c,\phi(c)]$ coplanar with $[a,b]$.

This is equivalent to $( \phi(a)- b | \phi(c) )=0$, which implies that $\phi(c)$ is contained in the hyperplane $\xi:=(\phi(a)- b)^\perp$ of $\Pu_k(\A)$. 

In order to have a unique possibility for $c$ (up to scalar products), it suffices that the subspace $\xi' = \phi(u^\perp)$ intersects $\xi \cap \phi(\Pu_k(\A))$ in a vector line. This is equivalent with $\xi' \neq \xi \cap \phi(\Pu_k(\A))$.

Assume by way of contradiction that we do have equality. Then $\phi(a)$, which is an element of $\xi'$, is also contained in $\xi$. This would imply that $( \phi(a)- b | \phi(a) ) = 0 $, which is equivalent to $(b, \phi(a) ) = | \phi(a) |^2$. As $|b| = |a|$, the Cauchy-Schwarz inequality now implies that $b  = l \phi(a)$, with $l \in k$ of unit norm. Plugging this back into $( \phi(a)- b | \phi(a) ) = 0 $ yields that $l=1$ or equivalently $b = \phi(a)$, which is a contradiction. This proves the proposition.
\qed

\begin{remark}
Another way to get the same conclusion is to interpret the condition $|a| = |b|$ for the lines through a fixed point as a Hermitian form on a $(\dim_k \A + \dim_k \B - 3)$-dimensional vector space over $k$, and obtain a generalized quadrangle in this way. 
\end{remark}

\subsection{Possibilities for $k$, $A$ and $B$}\label{sec:pos}
By Hurwitz's theorem on composition algebras we only have the following three possibilities for $k$, $\A$ and $\B$ (up to isomorphism).
\begin{align*}
k& = \R, \A= \HH, \B = \HH , \\
k& = \R, \A= \HH, \B = \OO, \\
k& = \C, \A= \OO, \B = \OO.
\end{align*}

In each of these cases one can construct the following compact flag-transitive automorphism group of the associated geometry $\Gamma$. Set $G = \Aut_k \A \times  \Aut_k \B$, and the action of an element $g \define (\alpha, \beta) \in G$ on the geometry to be 
\begin{itemize}
\item points: $ka \mapsto \alpha(ka)$
\item lines: $[a,b] \mapsto [\alpha(a),\beta(b)]$
\item planes: $\phi \mapsto \beta\phi\alpha^{-1}$.
\end{itemize}

It is straightforward to check this action is bijective and preserves incidence. As $\Aut_\R \HH \simeq \SO(3)$, $\Aut_\R \OO \simeq \mathsf{G}(2)$ and $\Aut_\C \OO \simeq \SU(3)$ (see \cite{CPP}), we obtain the following possibilities for $G$ (using the same ordering as before).

\begin{align*}
\SO(3) \times & \SO(3), \\
\SO(3) \times & \mathsf{G}_2, \\
\SU(3) \times & \SU(3). 
\end{align*}

\begin{remark}
The action is faithful except in the last case, for which the action has a kernel of size three.
\end{remark}

The group $G$ clearly admits a compact topology and acts continuously for the standard topology on $\Gamma$. Flag-transitivity is easily seen by the transitivity of $\Aut_k \A$ on $\Pu_k(\A)$, and the fact that, when assuming $\A$ to be a $k$-subalgebra of $\B$, each $k$-algebra injective morphism from $\A$ into $\B$ extends to an automorphism of $\B$.

We can hence apply the classification made in Theorem A of~\cite{Kra-Lyt:14} and conclude that the geometry $\Gamma$ is either covered by a building, or is one of two unique exceptional geometries of type $\mathsf{C}_3$. Our last two possibilities correspond with the two exceptional geometries, the geometry for the first possibility is therefore covered by a building of type $\mathsf{C}_3$.

\subsection{The case $G = \SO(3) \times SO(3)$}

For this case we provide an explicit description of the covering. The covering geometry is the projective quadric $\Delta$ defined by the following equation.

$$
Q(X_0,X_1,X_2,X_3,X_4,X_5,X_6)=X_0^2 + X_1^2 + X_2^2 - X_3^2 - X_4^2 -X_5^2 - X_6^2=0
$$

The group $\SO(4)$ acts on the last four coordinates, in particular on the unit sphere $\mathbb{S}^3$ in the three-dimensional space spanned by these. Let $H$ be the subgroup of left isoclinic rotations (which is isomorphic with $\SU(2)$ as well as the multiplicative group of unit quaternions, which themselves can be identified with the sphere $\mathbb{S}^3$). This group $H$ acts sharply transitive on $\mathbb{S}^3$.

Let $G$ be the group $\SO(3) \times \SO(4)$ (where $\SO(3)$ acts on the first three coordinates, and $\SO(4)$ is as before). Then $H\leq G$ and the quotient $G / H$ is isomorphic to $\SO(3) \times \SO(3)$.

The action of $G$ on $\Delta$ is well-understood in the context of Veronese and polar representations, in particular one has the following lemma.

\begin{lemma}\label{lem:transitiveG}
The group $G$ acts chamber-transitively on the projective quadric $\Delta$.
\end{lemma}

If we can show that the group $H$ does not fix any point, line or plane, then the quotient geometry $\Delta'$ is a $\mathsf{C}_3$ geometry on which the group quotient $G / H \simeq \SO(3) \times \SO(3)$ acts chamber-transitively.

\begin{lemma}
No non-identity element of $H$ maps a point of $\Delta$ to a collinear point on it.
\end{lemma}
\proof
Suppose that an element $g \in G$ maps some point $p$ of $\Delta$ to a collinear point $q$. By Lemma~\ref{lem:transitiveG} and the fact that $H$ is a normal subgroup of $G$ we may assume w.l.o.g. that the point $p$ is represented by $(1,0,0,1,0,0,0)$. The point $q$ can then be represented by $(1,0,0,a,b,c,d)$ where $a^2 +b^2 + c^2 + d^2 = 1$. 

In order for these points to be collinear on $\Delta$ we need that for the associated bilinear form $B(x,y)=\frac{1}{2}(Q(x+y)-Q(x)-Q(y))$ we have $B(p,q)=0$. 
This implies that $\frac{1}{2}(4 - (1+a)^2-b^2-c^2-d^2) = 0$. Combined with the previous condition $a^2 +b^2 + c^2 + d^2 = 1$ it follows that $a=1$ and hence $b=c=d=0$, so $p = q$. As $H$ acts sharply transitive on the sphere $\mathbb{S}^3$, we have that $g$ is the identity element. This proves the lemma.
\qed

\begin{prop}
The group $H$ acts freely on the projective quadric $\Delta$.
\end{prop}
\proof 
This follows directly from the previous lemma.
\qed

One can verify that the geometry and group action constructed in this section is the same one as in Section~\ref{sec:pos} by comparing the stabilizer of a chamber and the subsimplices of it.

\section{The full automorphism group of the geometry $\Gamma$}\label{section:auto}

This section is devoted to determining the full automorphism group $\Aut(\Gamma)$ of $\Gamma$, which will turn out to be very close to the compact group given in Section~\ref{sec:pos}.

\subsection{Statement of results}\label{section:fullclaim}
We claim that the full automorphism group is given by the following table (listed case-by-case):

\begin{align*}
k& = \R, \A= \HH, \B = \HH : \SO(3) \times \SO(3), \\
k& = \R, \A= \HH, \B = \OO : \SO(3) \times \mathsf{G}_2, \\
k& = \C, \A= \OO, \B = \OO : (( \SU(3) \times \SU(3)) / C_3) \rtimes C_2.
\end{align*}

The cyclic group $C_2$ consisting of two elements arises from complex conjugation. The cyclic group $C_3$ is the kernel of the map of the group $\SU(3)$ to its projectivization 
$\mathsf{PSU}(3)$, i.e. the center of $\SU(3)$, matrices $\zeta I$, where $\zeta$ is a third root of unity and $I$ is the identity matrix. 

From the results of Section~\ref{sec:pos} it easily follows that these are indeed automorphism groups of the associated geometries. It only remains to check if these are the full groups.

\subsection{Image into $\mathsf{P\Gamma L}(3, k)$}\label{section:image}

The set of points and the set of point shadows of lines (removing duplicity) form a projective plane $\Sigma$ over $k$, see also Remark~\ref{rem:flat}. 
Hence the automorphism group of $\Gamma$ can be mapped into the full automorphism group $\Aut(\Sigma)$ of this plane, which is $\mathsf{P\Gamma L}(3, k)$.
We aim to show that the image of $\Aut(\Gamma)$ into $\Aut(\Sigma)$ is contained in the centralizer of the polarity induced by the Hermitian inner product $(\cdot | \cdot )$, as expected.

As the compact automorphism group listed in Section~\ref{sec:pos} is flag-transitive, we may fix the point $x \define \<a \>$, and line $L \define [a',b]$, where $a, a', a''$ is an orthogonal $k$-basis for $\Pu_k(\A)$. 
Let $M$ be a second line with the same point shadow as $L$, which we may assume to be of the form $M \define [a', d]$.

Consider the lines $N$ through $x$ which are coplanar with both $L$ and $M$. 
Such a line is of the form $[c,f]$.

\begin{lemma}
The point shadows for varying $N$ are all identical if and only if $\< b \> = \<d \>$. In this case the common point shadow is formed exactly by the points corresponding to the vector lines in the $k$-span of $\{a,a'\}$.
\end{lemma}
\proof
The point shadows are identical if and only if there is a unique solution for $c$ up to scalar multiples.

The equations that have to be satisfied in order for a line $N \define [c,f]$ go through $x$ and being coplanar with both $L$ and $M$ are the following.
\begin{align*}
&c \in \< a', a'' \> \\
(a' | c)& = (b | f) = (d | f)
\end{align*}
Note that $c \notin \<a' \>$. 
If $c = a'l + a'' $ is a valid solution, then  $c = a'l - a'' $ is also. 
Hence in order to have unicity for $c$ up to scalar multiples we need that $l=0$, or equivalently that $(b |f)  =0$ for every $f$ such that $(b |f) = (d|f)$. 
The latter condition can be rewritten as $(b-d | f) =0$.
As the subspace $(b-d)^\perp$ is contained in $b^\perp$  if and only if the vector lines through $b$ and $d$ are identical, it follows that $\< b \> = \<d \>$.

The element $c$ is then a scalar multiple of $a''$, which determines the point shadow of $N$.
\qed

Note that the point shadow of such an $N$ forms a line through $x$ orthogonal to the point shadow of $L$. 
One can hence recognize the polarity $\rho$.

We conclude that the image of $\Aut(\Gamma)$ in $\Aut(\Sigma)$ is contained in $\mathsf{SO}(3)$ if $k = \R$, and $\mathsf{PSU}(3) \rtimes C_2$ if $k = \C$, where the cyclic group $C_2$ corresponds with complex conjugation.
Note that this corresponds with the possibilities listed in Section~\ref{section:fullclaim}.

\subsection{Kernel of the map}\label{section:kernel}

In this section we take a look at the kernel of the map into  $\mathsf{P\Gamma L}(3, k)$ as considered in Section~\ref{section:image}. 
If we can show that this is exactly $\Aut_k \B$, then the claim from Section~\ref{section:fullclaim} readily follows. 
As the kernels for the groups of automorphisms listed in Section~\ref{section:fullclaim} still act transitive on the set of planes of the geometry $\Gamma$, we may restrict ourselves to considering the stabilizer of one such plane. 
This choice of a plane corresponds with an embedding of $\A$ into $\B$, hence we may assume that $\A$ is a $k$-subalgebra of $\B$.

Let $H$ be the stabilizer of the plane corresponding with the (natural) embedding of $\A$ into $\B$ in the kernel of the map into  $\mathsf{P\Gamma L}(3, k)$.
So $H$ fixes every line of the form $[a,a]$ ($a \in \Pu_k\A$), as well as every point of $\Gamma$.

Let $a$ in $\Pu_k(\A)$ be an element of unit norm, and consider the set of lines $K \define \{[a,b] \vert b \in \Pu_k(\B), |a| = |b| =1 \}$, which is naturally parametrized by the elements in $\Pu_k(\B)$ of unit norm and stabilized by $H$. 
In the next few lemmas we investigate how $H$ acts on $K$, which via the parametrization corresponds with an action of $H$ on the elements of unit norm in $\Pu_k(\B)$.

\begin{lemma}\label{lemma:spaces}
The action of $H$ on the set $K$ preserves the intersection of vector subspaces of $\Pu_k(\B)$ (as a vector space over $k$) with the subset of elements of unit norm.
\end{lemma}
\proof
The lines $[a,b]$ in $K$ coplanar with some given line $[c,d]$ are determined by the linear equation $(a|c) = (b|d)$ over $k$ by Lemma~\ref{lem:coplanar}. We can therefore recognize the hyperplanes, and hence any vector subspace of $\Pu_k(\B) \cup \{0\}$ intersected with the elements of unit norm.
\qed

\begin{lemma}\label{lemma:ortho}
The action of $H$ on $K$ preserves orthogonality.
\end{lemma}
\proof
Consider a set $K' \define \{[c,d] \vert b \in \Pu_k(\B), |c| = |d| \}$ of lines with the same point shadow, where $a$ and $c$ are not scalar multiples of each other.

Then for a fixed line $[a,b] \in K$, the lines in $K$ which are coplanar with a line in $K'$ which is on its turn coplanar to $[a,b]$ are those lines $[a,f]$ where $f$ is such that there exists a $d \in \Pu_k(\B)$ with $|d| = |f| = |a|$ and $(a|c)=(b|d) = (f|d)$. 

If this is not the entirety of lines in $K$, it will certainly not contain those lines $[a,f]$ with $(c|f) = 0$. This characterizes exactly those lines.
\qed

\begin{lemma}\label{lemma:spread}
Let $a \in \Pu_k \A$ and $b \in \Pu_k \B$  such that $|a| = |b|$. It is then possible to reconstruct the lines of the form $[c,d]$ where $d$ is a scalar multiple of $b$ from the geometry.
\end{lemma}
\proof
We first consider the case where $(c|a) = 0$. 
The line $[a,b]$ is coplanar with a line $[c,f]$ ($f \in \Pu_k \B$ and $|f| = |c|$) if and only if $(f|d) =0$. 
As one can recognize orthogonality by Lemma~\ref{lemma:ortho}, it follows that one can reconstruct the lines $[c,d]$ where $d$ is a scalar multiple from $b$ from this.

The case where $(c|a) \neq 0$ follows from applying the previous case twice to an element $e$ in $\Pu_k \A$ such that $(e|a) = (e|c) =0$, which is always possible to find (as $(e|a) = (e|c) =0$ define two linear equations in a three-dimensional space).
\qed

\subsubsection{The case $\dim_k \B = 4$}\label{section:dim4}
We start by considering the case where $\B$ is four-dimensional over $k$, or, by the assumptions at the beginning of Section~\ref{section:kernel}, that $\A$ equals $\B$. We will show that $H$ acts trivially on any line and plane of the geometry, which proves the claim.

Assume that $[a, b]$ is line of the geometry $\Gamma$ (so $|a| = |b|$) not fixed by a certain element $h\in H$. 
By Lemma~\ref{lemma:spread} and the fact that $\A= \B$ the element $h$ can map this line only to a line $[a,bl]$, with $l \in k$ and $|l| = 1$. 

The lines of type $[c, c]$ ($c \in \Pu_k \A$) (which are fixed by $H$) coplanar with $[a,b]$ are those lines such that $(c|a) = (c|b)$. 
Note that there exists such $c$ up to scalar multiples, by considering the residue of a point incident with $[a,b]$, which is a generalized quadrangle. 
If $(c|a) \neq 0$, then we would have that $(c|b) = (c|b l) =(c|b) l$ implying that $l =1$ and that $[a,b]$ is fixed by $h$, which is a contradiction.

If there are only such coplanar lines $[c,c]$ where $(c|a) = 0$, then, by considering the residue of each point incident with $[a,b]$, the line $[a,b]$ is coplanar with each line of the form $[c,c]$ containing the point $\<a\>$.
So $a$ and $b$ are orthogonal to the same elements and are hence scalar multiples from each other.
Such a line is however also coplanar with fixed lines of the form $[c,d]$ where $c$ and $d$ are not scalar multiples of each other, and where $(c |a ) \neq 0$. 
As before one deduces that the line $[a,b]$ is fixed.

From this it follows that every line, and by extension every plane of $\Gamma$ will be fixed by $H$. This proves the claim.

\subsubsection{The case $\dim_k\B = 8$}
We now assume that $\dim_k\B = 8$. In particular this implies that $k = \R$, $\A = \HH$ and $\B = \OO$.  
The following lemma allows us to recognize quaternion subalgebras of $\B$.

\begin{lemma}\label{lemma:4sub}
The action of $H$ on $K$ preserves the intersection with four-dimensional $k$-subalgebras of $\B$.
\end{lemma}
\proof
Let $a'$ and $a''$ be elements of unit norm in $\Pu_k(\A)$ such that $(a|a') = 0$ and $a'' = a a'$. So $1, a, a'$ and $a''$ form an orthonormal $k$-basis for $\A$.

Let $1,b,b'$ and $b'' \define b b'$ be an orthonormal $k$-basis for a four-dimensional $k$-subalgebra of $\B$. 
We are going to construct the line $[a,b'' l'']$ (so $b''$ up to a scalar product) given the lines $[a,b]$ and $[a,b']$ in $K$. 
By Lemma~\ref{lemma:spread} we can construct lines $[a',bl]$ and $[a'',b'l']$ from these lines. The unique plane containing both $[a',bl]$ and $[a'',b'l']$ is given by the embedding defined by 
$$
\phi:    
\left\{
  \begin{array}{l}
    a \mapsto b'' ll' \\
    a'\mapsto bl \\
    a'' \mapsto b'l' \\
  \end{array}
\right.
$$
The unique line of the form $[a,f]$ incident with the plane is the line $[a, b''ll']$, as desired. The statement of the lemma now follows from Lemma~\ref{lemma:spaces}. \qed

This lemma allows us to consider the subgeometry of $\Gamma$ corresponding with the choice $\B = \A = \HH$, and conclude by the results of Section~\ref{section:dim4} that this subgeometry is completely fixed by $H$.

Recall the 3 cases we are considering and their proposed group of automorphism, see Section~\ref{section:fullclaim}.
\begin{align*}
k& = \R, \A= \HH, \B = \HH : \SO(3) \times \SO(3), \\
k& = \R, \A= \HH, \B = \OO : \SO(3) \times \mathsf{G}_2, \\
k& = \C, \A= \OO, \B = \OO : (( \SU(3) \times \SU(3)) / C_3) \rtimes C_2.
\end{align*}

The proposed group of automorphisms acts sharply transitively on the elements in $\B$ of unit norm and orthogonal to all elements in $\A$ (see for example~\cite[4.1]{Bae:02}), so fix such an element $b \in \B$.

Let $H'$ be the stabilizer of a line $[c,b]$ with $|c| = |b|$.
(Note that the choice of $c$ does not matter.)
Each element $a \in \Pu_k\A$ of unit norm generates a four-dimensional $k$-subalgebra of $\B$ together with $b$, which we can recover geometrically in the form of subgeometry $\Gamma'$ fixed by $H'$ by Lemma~\ref{lemma:4sub}. 
As we know that the group of automorphisms fixing all points of this subgeometry will be the group $\SO(3)$, and that $H$ additionally fixes the lines $[a,a]$ and $[a, b]$ of $\Gamma'$, we conclude that $H'$ fixes the subgeometry $\Gamma'$. 
Repeating this argument for different subgeometries we obtain that $H'$ acts trivially on $\Gamma$, whence the claim.

\section{On simple connectedness of the  geometry $\Gamma$}\label{section:simple}

In this section we prove that the geometry $\Gamma \define (V, \tau, *)$ of type $\mathsf{C}_3$ constructed in Section~\ref{sec:constr} is either simply connected, or covered by a building. 

\paragraph{Sketch of proof. --} We start by considering the edge-path group of $\Gamma$, which we use to study the universal cover $\widetilde{\Gamma}$. From Proposition~\ref{prop:control} (proved in Sections~\ref{section:red} up to~\ref{section:control}) we show that this cover admits a compact topology with a flag-transitive continuous group acting on it. Such a geometry is subject to the classification of Kramer and Lytchak (\cite{Kra-Lyt:14}) and is hence known, from which one concludes that the cover is either a building or the exceptional geometry itself.

\subsection{The edge-path group of $\Gamma$}
Instead of working with the fundamental group of $| \Gamma |$ it makes more sense to work with the equivalent edge-path group, which is more natural in the simplicial setting. We refer to~\cite[Ch. 44]{Sei-Thr:80} or to \cite{Spanier} for a detailed exposition. 

For us, an \emph{edge path} is a finite sequence of vertices in $V$ such that each two subsequent vertices are incident. (Note that such a sequence defines a unique sequence of edges on $| \Gamma |$.) The \emph{length} of a path is the number of vertices in it minus one.

We now introduce two kinds of \emph{elementary combinatorial deformations}. These are:

\begin{enumerate}[i.]
\item
Replacing a subpath $(u)$ with $(u,v,u)$, where $u$ and $v$ are incident, or the converse operation.
\end{enumerate}

\begin{enumerate}[i.]
\setcounter{enumi}{1}
\item
Replacing a subpath $(u,v)$ with $(u,w,v)$, where $w$ is incident with both $u$ and $v$ (which is equivalent with $\{u,v,w\}$ being a flag), or the converse operation.
\end{enumerate}

Two edge paths $\gamma$ and $\gamma'$ which can be transformed into each other by a finite number of elementary combinatorial deformations are said to be \emph{combinatorially homotopic}. An edge path is \emph{contractible} if it is combinatorially homotopic to a path of length zero. As combinatorial homotopies do not alter the begin or endpoint of a path, this implies that a contractible edge path begins and ends at the same point.

For a fixed vertex $v \in V$, the \emph{edge path group} $E(\Gamma, v)$ consists of the equivalence classes of combinatorially homotopic edge paths starting and ending at the vertex $v$. The group multiplication is defined as concatenation.

\subsection{The universal cover $\widetilde\Gamma$}

It will be useful to have a concrete model of the universal cover. We will do this by constructing a geometry $\widetilde\Gamma \define (\widetilde{V}, \widetilde{\tau}, \widetilde{*})$ using edge paths.

Fix a point $x \in V$. The set of vertices $\widetilde{V}$ consists of the equivalence classes of edge paths starting at $x$ under combinatorial homotopy. 
We set the type (i.e. its image under $\widetilde{\tau}$) of such an equivalence class of edge paths to be the type of the last vertex in any path in the class. 
Two equivalence classes of edge paths are incident if one can find an edge path of length $k$ in one class, and an edge path of length $k+1$ in the other class from which the first path can be obtained by removing the last vertex in the path.

Finally, the covering map $\rho: \widetilde{V} \to V$ is then defined by mapping an equivalence edge path to the last vertex of any path in it.

\subsection{Lifting the compact topology to the universal cover}

The goal of this section is to define a compact topology $\mathcal{T}$ on $\widetilde{V}$ starting from the compact topology on $V$. Our approach is based on Section 4 of~\cite{Fan-Gro-Tho2} and Section 8 of~\cite{Lyt:14}.

Proving the following proposition will be the subject of Sections~\ref{section:red} up to~\ref{section:control}.

\begin{prop}\label{prop:control}
For every  natural number $k$ there exists a natural number $D(k)$ with the following property. Any two edge paths $\gamma$ and $\gamma'$ in the universal cover $\widetilde{\Gamma}$ of length at most $k$ with the same extremities are homotopic  by a combinatorial homotopy consisting of at most $D(k)$ elementary combinatorial deformations.
\end{prop}

We fix a point $\tilde{x}$ in the universal cover $\widetilde\Gamma$, such that $x \define \rho(\tilde{x})$ is its image under the covering map $\rho$. 

We say that a sequence of vertices $(\tilde{v}_n)_{n \in \mathbb{N}}$ of vertices in $\widetilde\Gamma$ converges to a vertex $\tilde{v} \in \widetilde{V}$, if one can choose representative edge paths $(x, u_1^n, u_2^n, \dots, u_k^n)$ for each $\tilde{v}_n$, all of the same length $k$, and such that each of the sequences $(u_i^n)_{n \in \N}$ converges to some vertex $u_i$ (in the compact topology on $V$) and $(x, u_1, \dots, u_k)$ is a valid edge path representing $\tilde{v}$.

This notion of convergence defines closed sets and hence a topology $\mathcal{T}$ on $\widetilde{V}$. This topology is sequentially compact, as the set of paths of fixed length $k$ starting at $x$ is closed in $V^{k+1}$. The topology admits a dense countable subset as this also holds for $V$, and by the fact that every vertex in $\widetilde{V}$ can be represented by an edge path whose length is globally bounded.

\begin{remark}\label{secondcountable}
One can define the topology $\mathcal{T}$ in a slightly different way as follows. Consider the set of edge paths starting at a single point $p$ (denote this set by $\overline{V}$).  One can interpret this a subset of the countable product space $\Pi_{i=0}^\infty V$. On this set of edge paths we define equivalences as above, thus defining a quotient topology $\mathcal{T}'$. Since countable products of second countable spaces are again second countable, and since subspaces of second countable spaces are second countable and since the quotient map defined by this equivalence relation is open, $\mathcal{T}'$ is second countable. By Proposition 2.4 of \cite{Siwiec} a convergent sequence in the quotient topology lifts to a convergent sequence in the topology itself. As $\mathcal{T}$ was described by convergence and since first countability implies that a set is closed if and only if it is sequentially closed, we can conclude that $\mathcal{T}'$ is equal to $\mathcal{T}$. Hence $\mathcal{T}$ is second countable.
\end{remark}

This topology will have nice properties, as exhibited by the following lemmas.
\begin{lemma}\label{lem:onelimit}
Each sequence has at most one limit.
\end{lemma}
\proof
Let $(\tilde{v}_n)_{n \in \mathbb{N}}$ be a sequence of vertices, and let $(x, u_1^n, u_2^n, \dots, u_k^n)$ and $(x, {u'}_1^n, {u'}_2^n, \dots, {u'}_l^n)$ be two edge paths both representing $\tilde{v}_n$, one of length $k$, the other of length $l$, for every $n \in \mathbb{N}$. Suppose that the paths $(x, u_1^n, u_2^n, \dots, u_k^n)$ converge to an edge path $(x, u_1^n, u_2^n, \dots, u_k^n)$, and the paths $(x, {u'}_1^n, {u'}_2^n, \dots, {u'}_l^n)$ to an edge path $(x, u'_1, u'_2, \dots, u'_l)$. We have to proof that these two limit edge paths are combinatorially homotopic.

Each pair of paths $(x, u_1^n, u_2^n, \dots, u_k^n)$ and $(x, {u'}_1^n, {u'}_2^n, \dots, {u'}_l^n)$ are combinatorially homotopic by a bounded number of elementary combinatorial deformations (see Proposition~\ref{prop:control}). 
By sequential compactness these combinatorial homotopies (which can be thought of as sequences of edge paths) for each $n$ admit a convergent subsequence, which will be the desired combinatorial homotopy.
\qed

\begin{lemma}\label{lem:independent}
The topology $\mathcal{T}$ is independent of the choice of points $x$ and $\tilde{x}$.
\end{lemma}
\proof
For another point $\tilde{y} \in \widetilde{V}$, where $y \define \rho(\tilde{y})$, and any vertex $\tilde{v}$ in the cover, we can extend an edge path starting from $x$ and representing  $\tilde{v}$ to an edge path starting from $y$, by adding a fixed path from $y$ to $x$ in front. This operation does not influence the notion of convergence in $\mathcal{T}$, hence this topology is independent of the choice of $x$ and $\tilde{x}$.
\qed

\begin{prop}
The universal cover $\widetilde{\Gamma}$ of $\Gamma$ admits a compact metrizable topology with connected panels, invariant under a flag-transitive automorphism group.
\end{prop} 
\proof
The topology $\mathcal{T}$ is second countable by Remark \ref{secondcountable} and hence Hausdorff by Lemma~\ref{lem:onelimit}, is separable and sequentially compact, hence this topology is compact and metrizable. 
As the panels of $\Gamma$ are connected, the panels of $\Gamma'$ which are homeomorphic are also connected.

By the results of Section~\ref{sec:pos} we know that a compact flag-transitive group $G$ of continuous automorphisms acts on $\Gamma$. 
This group lifts to a group $\widetilde{G}$ acting flag-transitively on the universal cover $\widetilde{\Gamma}$, where the kernel of the map $\widetilde{G} \to G$ is the group of deck transformations of the universal cover. 
As the topology $\mathcal{T}$ is independent of the choice of $x$ and $\tilde{x}$ by Lemma~\ref{lem:independent}, the group $\widetilde{G}$ acts continuously.
\qed

We are now in the position to apply Theorem A of~\cite{Kra-Lyt:14} yielding that the geometry $\widetilde{\Gamma}$ is either a building of type $\mathsf{C}_3$, or is one of two possible exceptional geometries. 

So if 
$$
k = \R, A= \HH, B = \HH,
$$
then $\Gamma$ cannot be covered by one of the two exceptional possibilities (as the rank two residues do not match), hence it is not simply connected and covered by a building, and if
$$
k = \R, A= \HH, B = \OO
$$
or
$$
k = \C, A= \OO, B = \OO,
$$
so where $\Gamma$ is an exceptional geometry associated to the actions of $\SU(3) \times \SU(3)$ or $\SO(3) \times \mathsf{G}_2$ on the Cayley plane, then $\widetilde\Gamma$ would be homeomorphic to $\Gamma$ (as it is shown in~\cite{Kra-Lyt:14} that these cannot be covered by a building, and the fact that the rank two residues of the two cases are different). As the universal cover is simply connected by definition, we conclude that the geometry $\Gamma$ is simply connected and its own universal cover.

\subsection{Reducing planes from edge paths}\label{section:red}

The goal from now on is to produce a proof for Proposition~\ref{prop:control}. 
A first step is to reduce the set of edge paths one needs to consider.

If $(u,\pi,v)$ is an edge path where $u$ and $v$ are two vertices incident with a common plane $\pi$ of the geometry $\Gamma$, then this edge path is combinatorially homotopic to any edge path of the form $(u, w_1, w_2, w_3, \dots, w_{k-1}, v)$ completely contained in the residue of $\pi$ (so assuming that each of the $w_i$ is incident with $\pi$) by applying $k$ elementary combinatorial deformations.

Hence, if an edge path $\gamma$ of length $k$ does not start or end at a plane, we can find a combinatorially homotopic edge path of length at most $\floor{\frac{3k}{2}}$ containing only points and lines using at most $3\floor{\frac{k}{2}}$ elementary combinatorial deformations (as $\floor{\frac{k}{2}}$ is the maximum amount of planes in $\gamma$, and 3 is the diameter of a projective plane). 

Moreover if $(u_0, u_1, \dots, u_{k-1}, u_k)$ and $(v_0 \define u_0, v_1, \dots, v_{l-1}, v_l \define u_k)$ are two edge paths both completely contained in the residue of some plane $\pi,$, then they are combinatorially homotopic by applying at most $k+l$ elementary combinatorial deformations.

\subsection{Primitive edge paths}\label{section:prim}

We call an edge path of the form $(x,L,y,M,x)$, where $L$ and $M$ are two different lines through two different points $x$ and $y$, \emph{primitive}.

In this section we want to show that if a primitive edge path in $\Gamma$ is contractible, then it can be reduced to the trivial edge path $(x)$ by at most $K$ elementary combinatorial deformations, where $K$ is a universal constant. (Which is a special case of Proposition~\ref{prop:control}.)

Two points $x$ and $y$ of $\Gamma$ are said to be \emph{orthogonal} if the corresponding vector lines in $\Pu_k(\A)$ are orthogonal w.r.t. the inner product $(\cdot | \cdot)$.

\begin{lemma}
Every primitive edge path $(x, L, y, M, x)$ is homotopic (by using at most 12 elementary combinatorial deformations) to some primitive edge path $(x, L', y', M', x)$, where $x$ and $y'$ are orthogonal.
\end{lemma}
\proof
If $x$ and $y$ are orthogonal, there is nothing to prove, so assume that this is not the case. We start by noting that the polar line of $x$ (i.e. the points orthogonal to $x$), is a line of the projective plane obtained from projectivizing $\Pu_k(\A)$. 

Fix a line $N$ through $y$ which is coplanar to both $L$ and $M$. This is always possible as the residue of $y$ is a generalized quadrangle. The line $N$ intersects the polar line of $x$ in some point $y'$. Let $L'$, respectively $M'$, be the unique lines through $x$ and $y'$, in the unique plane containing $L$ and $N$, respectively $M$ and $N$. By Section~\ref{section:red} we know that $(x, L, y, M, x)$  and $(x, L', y', M', x)$ are combinatorially homotopic by performing twelve elementary combinatorial deformations. This proves the lemma. \qed

Note that this lemma implies that if $j$ elementary combinatorial deformations suffice for  a contractible edge path $(x, L, y, M, x)$, then $j+ 12$ elementary combinatorial deformations suffice for the homotopic edge path $(x, L', y', M', x)$. 

As the automorphism group $G$ of $\Gamma$ described in Section~\ref{sec:pos} acts transitively on pairs of orthogonal points of $X$, we may fix an orthogonal pair of points $x,y$, represented by vector lines $\<a \>$ and $\<a'\>$ with $a$ and $a'$ of unit norm in  $\Pu_k(\A)$ with $(a|a') =0$, and reduce the question to the problem whether there is a global bound for the needed number of elementary combinatorial deformations for a contractible primitive edge path $(x, L, y, M, x)$.

Let $a'' \define a a' \in \Pu_k(\A)$, then $a''$ is of norm one and orthogonal to both $a$ and $a'$, see Proposition 11.10 of \cite{CPP}. We now can represent $L$ and $M$ by respectively $[a'', b]$ and $[a'',c]$ (where $b$ and $c$ are of unit norm). By the action of the $k$-automorphisms of $\B$ the orbit of a primitive edge path $(x, L, y, M, x)$ under the corresponding automorphisms of $\Gamma$ (which fix $x$ and $y$) is completely determined by the inner product $( b \vert c )$. We call this inner product the \emph{PL-invariant} of $(x, L, y, M, x)$.

%The lines through both points are exactly those lines of the form $[(0,0,1), v]$ with $v \in \C^3$ and $\vert \vert v \vert\vert =1$. 
%
%The stabilizer of $x$ and $y$ in $\mathrm{Aut}(X)$ acts transitively on pairs of lines $([(0,0,1), v], [(0,0,1), w])$ where $\langle v,w \rangle$ is fixed. Assume that $L$ is coordinatized by $[(0,0,1), a]$ and $[(0,0,1), b]$. Clearly $\langle a, b \rangle \in [0,1[$.   

\begin{lemma}
If a primitive edge path $(x, L, y, M, x)$ with PL-invariant $l$ is contractible, where $l\in k \setminus \{ \pm 1\}$ and $| l |=1$, then there exists a contractible primitive edge path $(x, L', y, M', x)$ with PL-invariant $l'$ with $| l' | <1$. Note that in order for such an element $l$ to exist it is necessary that $k = \C$.
%
%we can find a point $y'$ and lines $L'$ and $M'$ such that $PL(x,y',L',M')=(0,lr\bar{r}+\bar{l}(1-r\bar{r})$ for any $r$ with $|r|\leq 1$.
\end{lemma}
\proof
If we represent $L$ by $[a'',b]$ (as before), then the condition on the PL-invariant implies that $M$ is represented by $[a'',bl]$.

Let $b'$ in $\Pu_k(\B)$ be of unit norm such that $(b| b') =0$.
The line $N \define [a, b']$ is incident with the point $y$ and coplanar with both $L$ and $M$ (for different planes), as
\[
( a' | a )  = ( a''| a ) = ( b | b') = ( bl | b') = 0,
\]
where we make use of Lemma~\ref{lem:coplanar}. 

If we put $b'' = b b'$, then $1, b, b'$ and $b''$ form an orthogonal $k$-basis of a subalgebra of $\B$. This leads to two embeddings

$$
\phi:    
\left\{
  \begin{array}{l}
    a \mapsto b' \\
    a'\mapsto b'' \\
    a'' \mapsto b \\
  \end{array}
\right.
$$
and
$$
\psi:    
\left\{
  \begin{array}{l}
    a \mapsto b' \\
    a' \mapsto   b''\bar{l} \\
    a'' \mapsto  bl \\
  \end{array}
\right.
$$
of $\A$ into $\B$ (similar as in the proof of Proposition~\ref{prop:map}). The corresponding planes contain respectively the lines $L$, $N$ and $M, N$. Set $d \define \frac{1}{\sqrt{2}} (a'+a'')$ (note that this is of norm one and orthogonal to $a$). The lines $L' \define [d, \phi(d)] = [d, \frac{1}{\sqrt{2}}  (b  + b'')]$ and $M' \define [d, \psi(d)] = [d, \frac{1}{\sqrt{2}}  ( bl  +  b'' \bar{l})]$ lie in a common plane with respectively $L, N$ and $M, N$, contain $x$ and a common point $y'$ of $N$, which is orthogonal to $x$. By this reasoning $(x, L, y, M, x)$ will be combinatorially homotopic with $(x, L', y', M', x)$, in particular the latter will be contractible. 

By using a suitable automorphism of our geometry we can map this last primitive pl-path to a pl-path $(x, L', y', M', x)$, with PL-invariant $(\frac{1}{\sqrt{2}}  (b  + b'') \vert  \frac{1}{\sqrt{2}}  ( bl  +  b'' \bar{l}) ) = \frac{1}{2} (l + \bar{l} ) = \operatorname{Re} (l) < 1$. This proves the lemma.
\qed

We can distinguish between two possibilities at the moment, either there exists some contractible primitive edge path $(x, L, y, M, x)$ (where $L \neq M$) with PL-invariant different from $-1$, or every such contractible edge path has PL-invariant $-1$. (If the PL-invariant would be $1$, then we would have that $L = M$.)

In the second case the group of automorphisms fixing $x$ and $y$ acts transitively on such paths, hence there exists a global bound on the number of elementary combinatorial deformations needed to reduce a contractible primitive edge path to the trivial path in this case. 

In the remainder of the section we handle the first case.

\begin{lemma} \label{lem:diam}
Fix an $l \in k$. Let $\Gamma$ be the directed graph where the vertices are  the norm one elements in $\Pu_k (\B)$, and where $(b,c)$ is a directed edge if $( b |c ) = l$. If $\vert l \vert < 1$, then $\Gamma$ is connected and of finite diameter.
\end{lemma}
\proof
The proof will go in different steps, each one asserting that if $(b \vert c )$ equals some $l'$, then $c$ can be reached from $b$ using a finite number of directed edges. We call such an $l'$ \emph{valid}. By the action of the automorphism group of $\B$ over $k$, the number of edges needed depends only $l'$, not on $c$.

First we set $l' = l^2$. Then we may assume without loss of generality, that 
$$
c \define  b l^2 +  b' r
$$
where $b$, $b'$ and $b''$ are mutually orthogonal and of norm 1, and such that $| r | = \sqrt{1 - | l | ^2}$.  Set $d \define bl + b '' r$, then $(b,d)$ and $(d,c)$ are directed edges. We conclude that $l'$ is valid.

By repeating this argument we obtain that any power $l^n$, where $n$ is a power of two, is valid. Note that $| l^n |$ approaches zero.

We now claim that $0$ is valid. To see this we can take elements 
\begin{align*}
b& \\
d& \define b l^n  + b' {\sqrt{1- | l^n |^2} } b' \\
c& \define b' \frac{l^n}{\sqrt{1- | l^n |^2} } b' + b'' s  
\end{align*}
where $b$, $b'$ and $b''$ are as before, and $s$ is such that $|c| = 1$. The latter is only possible when $ | \frac{l^n}{\sqrt{1- | l^n |^2} } | \leq 1$, which is true for a large enough $n$ (which we fix for the remainder of the proof). This proves that $0$ is valid. In particular if $( b \vert c )$ = 0, then $c$ can be reached by at most $2n$ steps from $b$ (as $(b|d) = l^n$ and $(d|c) = l^n$).

This implies that the directed graph $\Gamma$ is connected and has diameter at most $4n$, as one can find for every pair of elements a mutually orthogonal element.\qed

We now make use of the following observation. If the primitive edge paths $(x, L, y, M, x)$ and $(x, M, y, N, x)$, both containing $x$ and $y$, are both contractible, then $(x, L, y, N, x)$ will be contractible as well. Moreover if the needed combinatorial homotopies of $(x, L, y, M, x)$ and $(x, M, y, N, x)$ both consist of at most $j$ elementary combinatorial deformations, then   there exists a combinatorial homotopy from  $(x, L, y, N, x)$ to the constant path $(x)$ consisting of at most $2j+2$ elementary combinatorial deformations. (There are two of them needed to obtain $(x, L, y , M, x, M,y, N, x)$, and then $2j$ deformations to reduce the subpaths.)

Combining this observation with Lemma~\ref{lem:diam} implies that every primitive edge path is contractible using at most $K$ elementary combinatorial deformations (where $K$ is a constant). 

Our claim holds in both cases, and hence in general.

\begin{remark}
In the first case one can, at this point, directly show that $\Gamma$ is simply connected. This is done by exploiting the fact that every primitive edge path is contractible.
\end{remark}

\subsection{Homotopy control}\label{section:control}

Throughout this section we work with edge paths in the universal cover $\widetilde{\Gamma}$ of $\Gamma$, in particular any two edge paths with same begin and endpoint are homotopic.

We begin with some lemmas on edge paths consisting only of points and lines.

\begin{lemma}\label{lemma:pinching}
Let $\gamma \define (x, L, y, M,z)$ be some edge path, and $L'$ a line (where $L \neq L'$) through $x$ coplanar with $L$. Then there exists an edge path $\gamma' \define (x,L',y',M',z)$ which can be obtained from $\gamma$ by 12 elementary combinatorial deformations.\end{lemma}
\proof
Let $\pi$ be the plane incident with both $L$ and $L'$. 
There exists some line $N$ incident with $\pi$ and a plane $\xi$ incident with both $N$ and $z$ (by Lemma~\ref{lem:c3}). 
Let $y'$ be the intersection point of $L'$ and $N$ and let $M'$ be the line incident with $\xi$ and both  points $y'$ and $z$.

With six elementary combinatorial deformations we get from $(x, L, y, M,z)$ to $(x, L',y',N,y, M,z)$, and with another six to $(y, L', y', M',z)$.
\qed

\begin{lemma}\label{lem:1and2}
Let $(x, L, y, M,z)$ and $(x, N, z)$ be two edge paths starting and ending at the same point, one of length four, the other of length two. Then we can obtain one path out of the other with at most $K+30$ elementary combinatorial deformations.
\end{lemma}
\proof
By applying Lemma~\ref{lemma:pinching} at most twice, we may assume that $L$ and $N$ are identical (at the cost of at most 24 elementary combinatorial transformations). As the path $(x, N, z)$ is combinatorially homotopic with $(x, N, y, N, z)$ (using six elementary combinatorial deformations), we are reduced with the question of transforming the subpath $(y,M,z)$ to the subpath $(y,N,z)$, which takes at most $K$ elementary combinatorial homotopies.
\qed

\begin{lemma}\label{lem:2and2}
Let $ (x, L, y, M,z)$ and $(x, L', y', M', z)$ be two edge paths starting and ending at the same point, both of length four. Then we can obtain one path out of the other with at most $K+55$ elementary combinatorial deformations.
\end{lemma}
\proof
As before we may assume by Lemma~\ref{lemma:pinching}  that $L$ and $L'$ are identical (at the cost of at most 24 elementary combinatorial deformations). One additional elementary combinatorial deformation transforms the path $(x,L,y,M,z)$ into $(x,L,y', L, y,M,z)$, which reduces the problem to the subpaths $ (y', L, y,M,z)$ and $(y', M', z)$, which is covered by Lemma~\ref{lem:1and2}. \qed

\begin{lemma}\label{lem:3}
If $\gamma \define (x_0, L_1, x_1, \dots , x_l)$ is an edge path (containing only of points and lines) of length $2l$, with $l \geq 3$ then there is an edge path $\gamma'$  of length $2l-2$ which one can be obtained from $\gamma$ by applying at most $K+56$ elementary combinatorial deformations.
\end{lemma}
\proof
Let $y$ be a point on $L_3$ collinear to $x_0$ (which exists by Lemma~\ref{lem:c3}), and $M$ a line incident with both $y$ and $x_0$.  The subpath $(x_0, L_1, x_1, L_2, x_2)$ is now combinatorially homotopic to the path $(x_0, M, y, L_3, x_2)$, using at most $K+55$ elementary combinatorial deformations, by Lemma~\ref{lem:2and2}. 

So the edge path $(x_0, L_1, x_1, L_2, x_2, L_3,  x_3, \dots ,x_l)$ is combinatorially homotopic to $(x_0, M, y, L_3, x_2, L_3,  x_3, \dots ,x_l)$, which in turn is combinatorially homotopic to the edge path $(x_0, M, y, L_3, L_3,  x_3, \dots ,x_l)$ of length $l-1$. \qed

The following lemma yields a more general homotopy control.

\begin{lemma}\label{lem:contr1}
For every  natural number $k$ there is a number $C(k)$ with the following property. Any two edge paths $\gamma$ and $\gamma'$ starting at the same point and ending at the same point, and containing only points and lines, in the universal cover $\widetilde{X}$ of length at most $k$ are homotopic  by a homotopy consisting of at most $C(k)$ elementary combinatorial deformations.
\end{lemma}
\proof
Note that the length of such paths is necessarily even.

For $k = 0$, $C(0) = 0$ trivially suffices. For $k=2$ we need $C(2) = K$, and for $k=4$ we can set $C(4) = K+55$ by Lemmas~\ref{lem:1and2} and~\ref{lem:2and2}.

For $k \geq 6$, we can reduce edge paths of length at most $k$ to paths of length at most four by applying at most $\frac{(k-4)(K+56)}{2}$ elementary homotopies by repeated application of Lemma~\ref{lem:3}, which leaves us in the previously handled case of $k=4$. It hence suffices to set $C(k)$ to be $(k-4)(K+56) + K+55$ in this case.
\qed

Finally, the next proposition removes the conditions concerning types, arriving at our claim in Proposition~\ref{prop:control}.

\begin{prop}
For every  natural number $k$ there is a number $D(k)$ with the following property. Any two edge paths $\gamma$ and $\gamma'$ with the same extremities are homotopic  by a homotopy consisting of at most $D(k)$ elementary combinatorial deformations.
\end{prop}
\proof
Let $u$ and $v$ be respectively the begin and end vertex of the edge paths $\gamma$ and $\gamma'$. 
If $u$ is a line or a plane, then we can apply an elementary combinatorial deformation to the beginning subpath $(u)$ of both $\gamma$ and $\gamma'$ to the subpath $(u,x,u)$ where $x$ is a point incident with $u$. 
If we apply a similar operation to the end vertex $v$, we can by omitting the extremities reduce the question to two edge paths $\gamma$ and $\gamma'$ of length at $k+2$ starting and ending at a point (at the cost of at most four elementary combinatorial deformations).

By applying Section~\ref{section:red} we may assume that $\gamma$ and $\gamma'$ are edge paths of length at most  $\floor{\frac{3k+6}{2}}$  consist of only points and lines (using at most $3\floor{\frac{k+2}{2}}$ elementary combinatorial deformations).

At this point we can apply Lemma~\ref{lem:contr1} and conclude that $D(k) = C(\floor{\frac{3k}{2}}) + 4 + 6\floor{\frac{k+2}{2}}$ is sufficient. 
\qed

\renewcommand{\refname}{Bibliography}
\bibliography{literaturliste}
\bibliographystyle{alpha}

\end{document}